\documentclass[reqno]{amsart}
\usepackage{amsmath}
\usepackage{amscd}
\usepackage{graphics}
\usepackage{latexsym}
\usepackage{color}
\usepackage{verbatim}
\usepackage{extarrows}
\usepackage{latexsym}
\usepackage{enumerate}
\usepackage[margin=1in]{geometry}
\usepackage{MnSymbol} 
\usepackage{graphicx}
\usepackage{tikz}
\usepackage[normalem]{ulem}
\usepackage[balance,small,nohug]{diagrams}
\usepackage[all]{xy}
\theoremstyle{plain}
\newtheorem{theorem}{Theorem}[section]

\newtheorem{corollary}[theorem]{Corollary}

\theoremstyle{definition}

\newtheorem*{thm1}{Theorem A}

\theoremstyle{remark}

\newcommand{\PP}{\mathbb{P}}

\newcommand{\PIC}{\mathrm{Pic}}

\def\PP{{\mathbb P}}

\begin{document}

\title{On the position of nodes of plane curves} 
\author{C\'esar Lozano Huerta}
\author{Tim Ryan}

\address{Universidad Nacional Aut\'onoma de M\'exico\\
Instituto de Matem\'aticas \\
Oaxaca, Mex.}
\email{lozano@im.unam.mx} 

\address{University of Michigan\\
Department of mathematics \\
Ann-Arbor, USA.}
\email{rtimothy@umich.edu} 


\maketitle
\begin{abstract}
The Severi variety $V_{d,n}$ of plane curves of a given degree $d$ and exactly $n$ nodes admits a map to the Hilbert scheme $\PP^{2[n]}$ of zero-dimensional subschemes of $\PP^2$ of degree $n$. This map assigns to every curve $C\in V_{d,n}$ its nodes. For some $n$, we consider the image under this map of many known divisors of the Severi variety and its partial compactification. We compute the divisor classes of such images in $\PIC(\PP^{2[n]})$ and provide enumerative numbers of nodal curves. We also answer directly a question of Diaz-Harris about whether the canonical class of the Severi variety is effective. 
\end{abstract}

\medskip
\section{Introduction}
\noindent
Plane curves are among the most classical objects in algebraic geometry. There is a solid understanding of them; individually as well as of many of their families. The family $V_{d,n}$ parameterizing irreducible nodal plane curves of degree $d$ with precisely $n$ nodes is particularly interesting. 
This family has been extensively studied and many of its basic properties are known \cite{DH,DH2,Ser}. A little less explored aspect of this family $V_{d,n}$ is: what is the geometry described by the nodes of the curves in it? This note studies this question.

\medskip\noindent
In many cases, the nodes of the curves in the family $V_{d,n}$ may simply be in general position. This often occurs each time the dimension of $V_{d,n}$ is at least $2n$, which is the dimension of the family of $n$ points on the plane \cite{AC}. In other words, the rational map $f:V_{d,n}\dashedrightarrow \PP^{2[n]}$ that assigns to each curve $C\in V_{d,n}$ its nodes is dominant. Here $\PP^{2[n]}$ stands for the Hilbert scheme of $n$ points on $\PP^2$. Even in this case,  the following is unknown as far as we are aware, what are the loci in $\PP^{2[n]}$ that occur as the nodes of families of irreducible curves? This note answers many cases of this question. Let us describe the situation precisely.

\medskip\noindent
We focus on a partial compactification of $V_{d,n}$, denoted by $W$, studied in \cite{DH, DH2}. Following these two papers, we will call $W$ the Severi variety. This space is smooth and one can describe readily many divisors in it. 
Furthermore, if $D\subset W$ is an effective divisor, then $f(D)$
might be a divisor as well. We will compute the divisor classes in $\PIC(\PP^{2[n]})$ of such images in the cases when the map $f:W\dashedrightarrow \PP^{2[n]}$ is birational. This will tell us the classes of divisors along the closure of the family of $n$ points in general position that occur as nodes of irreducible curves. Let us define the divisors in $W$ whose image we will compute:
\begin{itemize}
    \item the locus $CP$ of reduced and irreducible curves of genus $g
    $ with $n$ nodes containing a fixed point $p$;
    \item the locus $NL$ of reduced and irreducible curves of genus $g
    $ with $n$ nodes and one node located on a fixed line $L$;
    \item the locus $TN$ of reduced and irreducible curves of genus $g
    $ with $n-2$ nodes and one tacnode;
    \item the locus $TR$ of reduced and irreducible curves of genus $g
    $ with $n-3$ nodes and an ordinary triple point;
    \item the locus $\Delta_{0,1}$ of reduced curves of geometric genus $g-1$ with $n+1$ nodes, having two irreducible components, one of which maps to a line. 
\end{itemize}

\medskip\noindent
All of the classes computed in \cite{DH} are contained in a subspace of the Picard group $\PIC(W)_{\mathbb{Q}}$ generated by the classes of these $5$ divisors. In order to compute classes of the images of all divisors from \cite{DH},
it then suffices to compute the classes of the images of these $5$ generators. This is what we will do.

\medskip\noindent
In order to state the theorem, we recall the standard basis of the Picard group of $\mathbb{P}^{2[n]}$.
The class $H[n]\in \PIC(\PP^{2[n]})$ represents the subvariety of subschemes of length $n$ in $\PP^2$ whose support intersects a fixed line. 
Similarly, $B[n]$ is the class of the family of non-reduced subschemes of length $n$. The proof of the following result is presented in the next section. 

\begin{thm1}
\label{BASIS}
Let $W$ be the Severi variety of curves of degree $d$ and $n$ nodes. Suppose that the forgetful map $f:W\dashedrightarrow \PP^{2[n]}$, which assigns to a curve its nodes, is birational. Then the images under $f$ of $TN$, $TR$, $CP$, $NL$ and $\Delta_{0,1}$ have the following classes in $\PIC(\PP^{2[n]})$: 
\begin{equation}
    \begin{aligned}
f_*(TN) &= B[n],\\ f_*(TR) &=  0,\\ 
f_*(CP) &= (3d-3)H[n]-\tfrac{5}{2}B[n],\\ f_*(NL) &= H[n]\text{, and }\\ f_*(\Delta_{0,1}) &= 0.
    \end{aligned}
\end{equation}
\end{thm1}

\medskip\noindent
As a corollary, we will compute the classes in $\PIC(\PP^{2[n]})$ that come from the following subvarieties in the Severi variety $W$. We refer the reader to \cite{DH} for a detailed exposition about them.

\begin{itemize}
    \item the locus $CU$ of reduced and irreducible curves of genus $g
    $ with $n-1$ nodes and one cusp;
    \item the locus $TL$ of curves tangent to a fixed general line;
    \item the locus $FP$ of curves with a flex line passing through a fixed general point; 
    \item the locus $FL$ of curves with a flex located somewhere on a fixed general line;
    \item the locus $FN$ of curves with a flecnode;
    \item the locus $NP$ of curves such that the tangent line to a branch of a node passes through a fixed general point;
    \item the locus $HF$ of curves with a hyperflex;
    \item the branch locus $\text{BR}_N$ of the divisor $N$, in the universal family of $W$, whose general point consists of a curve and one of its nodes;
    \item the branch locus $\text{BR}_T$ of the divisor $T$, in the universal family of $W$, whose general point consists of a curve and a point on that curve such that the tangent line at that point passes through a fixed general point;
    \item the branch locus $\text{BR}_F$ of the divisor $F$, in the universal family of $W$, whose general point consists of a curve and one of its flexes;
    \item the canonical divisor $K_W$.
\end{itemize}

\medskip\noindent
The computations of \cite{DH} along with Theorem A yield the classes in $\PIC(\PP^{2[n]})$ of the images of the previous divisors. We defer the proof to the next section.  

\medskip\noindent
\textbf{Notation:} 
For simplicity in what follows, if $D\in \PIC(W)$ is a divisor class, then we will denote its image in the Picard group of the Hilbert scheme $\PIC(\PP^{2[n]})$ by $D$ instead of writing $f_*(D)$ all the time.


\begin{corollary}\label{CLASSES}
Let $W$ be the Severi variety and the forgetful map $f$ as in Theorem A. We have the following classes in $\PIC(\PP^{2[n]})$: 
\begin{itemize}
    \item CU = $(d^3+2d^2-d-6)H[n]-(\frac{5}{6}d^2+\frac{5}{2}d+2)B[n]$;\\
    \item TL = $2(3d^2-6d+2)H[n]-5(d-1)B[n]$;\\
    \item FP = $(\frac{15}{2}d^3-30d^2+\frac{39}{2}d+6)H[n]-\frac{25}{4}d(d-3)B[n]$;\\
    \item FL = $6(3d^2-6d+2)H[n]-15(d-1)B[n]$;\\
    \item FN = $(\frac{5}{2}d^3+5d^2-\frac{5}{2}d-18)H[n]-(\frac{25}{12}d^2+\frac{25}{4}d+2)B[n]$,\\
    \item NP = $\frac{1}{2}(d^3+2d^2-d-4)H[n]-\frac{5}{12}d(d+3)B[n]$;\\
    \item HF = $4 (11 d^3 - 68 d^2 + 79 d - 9)H[n]-\frac{2}{3}(55 d^2 - 285 d + 132)B[n]$;\\
    \item $\text{BR}_N$ = $(d^3+2d^2-d-6)H[n]-\frac{5}{6}d(d+3)B[n]$,\\
  \item $\text{BR}_T$ = $(13d^3-64d^2+53d+12)H[n]-\frac{1}{6}(65 d^2 - 255 d + 36)B[n]$,\\
    \item $\text{BR}_F$ = $6 (13 d^3 - 79 d^2 + 83 d + 2)H[n]-(65d^2-330d+111)B[n]$\text{, and}\\ 
    \item $\text{K}_W$ = $-3H[n]$.\\
\end{itemize}
\end{corollary}

\medskip\noindent
Diaz and Harris computed the canonical class of $W$:
$$K_W=-\tfrac{3}{5}A+\tfrac{3}{5}B+\tfrac{11}{12}C-\tfrac{13}{12}\Delta,$$ where the classes $A,B,C,\Delta$ are defined in the next section following \cite{DH}. 
In our present circumstances, Corollary \ref{CLASSES} asserts that the image of $K_W$ in $\PIC(\PP^{2[n]})$ is equal to $$f_*(K_W)=-3H[n]=K_{\PP^{2[n]}}.$$ Since this class fails to be pseudo-effective, it follows that $K_W$ cannot be effective. 
\begin{corollary}
Let $W$ be the Severi variety with values as in Theorem A. Then, $K_W$ is not effective. 
\end{corollary}

\noindent
This corollary answers a question in \cite[pag. 10]{DH} for the cases of this note.

\medskip\noindent
We finish this note in Section \ref{3} with an application to enumerative geometry. We have pushed divisors to the Hilbert scheme, thus we can now perform intersection theory on it.  We will intersect moving curves on $\PP^{2[n]}$ with the divisor classes we have computed. These intersection numbers will provide enumerative information about the position of the nodes of irreducible curves.

\section{Computations of the divisor classes} 
\medskip\noindent
In this section we compute the classes in Theorem A and in Corollary \ref{CLASSES}. Let $(d,n)$ be integers such that the forgetful map $f:W\dashedrightarrow \PP^{2[n]}$ is birational throughout. This is the if we have that $6n=d^2+3d$, except $(d,n)\neq (6,9)$, due to \cite{Tre,AC}.
\noindent
\begin{proof}[Proof of Theorem A]
We next note that $CP$ is precisely one of the Severi divisors defined in \cite{LR}.
In that paper, the class is computed to be $CP = (3d-3)H[n]-\frac{5}{2}B[n]$.

Straightforwardly, the closure of the image of $\Delta_{0,1}$ is the locus of collections of points with at least $d-1$ collinear points.
Since $d\geq 5$, this locus is not divisorial.
Thus, $\Delta_{0,1}$ is contracted, and the class of the image is 0.

Similarly, the closure of the image of $TR$ is the locus of collections of points supported on at most $n-2$ points. 
Since this locus is not divisorial in $\PP^{2[n]}$, $TR$ is contracted, and the class of its image is 0.

The remaining two divisors do map to divisors and have not previously been computed.
Treger \cite{Tre} proved that the map from the Severi variety to the Hilbert scheme is birational into its image in our present circumstances.
In fact, the map restricted to each of these two divisors is again birational.
Indeed, the map $f$ is regular along the divisors $NL$ and $TN$, cf. \cite[pag. 3]{DH} and the Severi variety $W$ is smooth along these two divisors. Since containing nodes at general points or on a general line are independent linear conditions inside the projective space of degree $d$ curves (similarly, containing nodes at general points or having a tacnode at a general point are independent linear conditions inside the projective space of degree $d$ curves), then the birational map $f$ over $V_{d,n}$ extends to a birational map over $NL$ and $TN$.
Thus, the pushforward class of each divisor is precisely the class of the reduced structure on the closure of its image.

The closure of the image of $TN$ is precisely the locus of collections of points supported on at most $n-1$ points.
This closure is precisely the exceptional divisor of the Hilbert-Chow map, and we have $TN = B[n]$.

Finally, the closure of the image of $NL$ is the locus of collections of points with at least one of the points on a fixed line.
Thus, $NL = H[n]$.
\end{proof}

\begin{proof}[Proof of Cor. 1.1] 
The map $f$ is proper between schemes and birational into its image \cite{Tre}. Moreover, it induces a morphism between Picard groups,  $$f_*:\PIC(W)_{\mathbb{Q}}\rightarrow  \PIC(\PP^{2[n]})_{\mathbb{Q}}.$$
Indeed, there is an injective morphism from the Picard group to the Chow group $\pi:\PIC_\mathbb{Q}(W)\rightarrow A_{2n-1}(W)\otimes \mathbb{Q}$ as the space $W$ is smooth \cite{DH}. The map $\tilde{f}_*:A_{2n-1}(W)\rightarrow  A_{2n-1}(\PP^{2[n]})$ is the morphism on Chow groups. Since there is an isomorphism $A_{2n-1}(\PP^{2[n]})\otimes \mathbb{Q}\xrightarrow{\cong} \PIC_{\mathbb{Q}}(\PP^{2[n]})$, then the claim follows by taking the composition $\pi\circ\tilde{f}_*$ . Consequently, $f_*(\alpha Q+Q')=\alpha f_*(Q)+f_*(Q')$ for any $\alpha\in \mathbb{Q}$ and classes $Q,Q'$. \cite[Section 1.4]{Fu98}.

In \cite{DH}, they compute the classes of each of the divisors from the statement of the corollary as linear combinations of the divisors $A$, $B$, $C$, $\Delta$, and $\Delta_{0,1}$.
We have already defined $\Delta_{0,1}$; 
let us define the rest now.
The class $\Delta$ represents the locus of reduced curves of geometric genus $g-1$, having at most two
irreducible components, with $n + 1$ nodes.
The remaining three of these are defined using divisors from the universal family $\mathcal{C}$ over the Severi variety.
The universal family comes equipped with two projections $\eta: \mathcal{C}\to \mathbb{P}^2$ and $\pi: \mathcal{C} \to W$.
If we let $\omega$ be the first Chern class of the relative dualizing class of $\mathcal{C}$ over $W$ and $D = \eta^*(c_1\left(\mathcal{O}_{\mathbb{P}^2}(1)\right))$, then we can define 
$A = \pi_*\left(D^2 \right)$, $B = \pi_*\left(D\cdot \omega \right)$, and $C = \pi_*\left(\omega^2 \right)$.

Using Mathematica and Macaulay2, the divisor classes of the corollary, written in the basis $A$, $B$, $C$, $\Delta$, and $\Delta_{0,1}$, can be converted to the basis of $TN$, $TR$, $CP$, $NL$, and $\Delta_{0,1}$. 
The change of basis is achieved by the following matrix \[\begin{bmatrix}
0&0&2&4&0\\0&0&9&15&0\\1&2d-3&
\frac{1}{6}(13d - 9)(d - 6)
&\frac{1}{6}(11d^2-57d+18)&0\\0&-2&-5d+18&-7d+18&0\\0&0&0&0&1\\
\end{bmatrix}.\]
The result follows from the linearity of the pushforward map and Theorem A.
\end{proof}

\section{Application to enumerative geometry}\label{3}
\noindent 
We may apply the results of Theorem A and Corollary \ref{CLASSES} to answer enumerative questions. To this end, we will consider the classes of closures of images $f(D)$ of divisors $D$ on the Severi variety $W$. The intersection of such classes with curves may occur along the closure; outside the image $f(D)$. In order to avoid this situation, we will consider moving curves. In other words, curve classes in the Mori cone $\overline{\mathrm{NE}}(\PP^{2[n]})$ whose representatives cover an open dense subset of $\PP^{2[n]}$. 
We know a geometric description of some such curves due to \cite{Hu14, BDPP}. 

\medskip\noindent
Let us define two moving curves in $\PP^{2[n]}$. Let $C_1$ be a curve in $\PP^{2[n]}$ defined as the collection of points containing $n-1$ general fixed points and whose final point varies on a general fixed line. In order to define the second moving curve, let us first write $n=\tfrac{r(r+1)}{2}+s$, with $0\leq s \leq r$. With this notation, a general $\Gamma\in \PP^{2[n]}$ lies on a smooth curve $C$ of degree $r$. Let us consider the curve $C_2$ induced in $\PP^{2[n]}$ by moving $\Gamma$ in a general pencil in the linear system $|\mathcal{O}_C(\Gamma)|$. In fact, if $s/r\in \Phi=\{\frac{s}{r}:(1+\sqrt{5})s > 2r\} \cup \{\frac{0}{1},\frac{1}{2},\frac{3}{5},\frac{8}{13},\cdots\}$, 
then \cite{Hu14} shows that these curves, $C_1$ and $C_2$, generate the moving cone of curves. 


\begin{table}
    \begin{tabular}{|c|c|c|}\hline
          & $C_1$ & $C_2$ \\ \hline  
          $TN$ & 0 & $(r-1)(r-2)-2+2n$\\\hline 
          $TR$ & 0 & 0\\\hline
          $CP$ & $3d-3$ & $(3d-3)r-\tfrac{5}{2}((r-1)(r-2)-2+2n)$ \\\hline
          $NL$ & 1 & $r$\\\hline
          $\Delta_{0,1}$ & 0&0\\\hline
          $CU$ &$d^3+2d^2-d-6$& 
          $(d^3+2d^2-d-6)r-(\frac{5}{6}d^2+\frac{5}{2}d+2)((r-1)(r-2)-2+2n)$
 \\\hline
          $TL$ & $2(3d^2-6d+2)$ & $2(3d^2-6d+2)r-5(d-1)((r-1)(r-2)-2+2n)$\\\hline
          $FP$ & $\frac{3}{2}(5 d^3 - 20 d^2 + 13 d + 4)$ & $(\frac{15}{2}d^3-30d^2+\frac{39}{2}d+6)r-\frac{25}{4}d(d-3)((r-1)(r-2)-2+2n)$\\\hline
          $FL$ & $6(3d^2-6d+2)$ & $6(3d^2-6d+2)r-15(d-1)((r-1)(r-2)-2+2n)$\\\hline
          $FN$ & $\frac{1}{2}(5 d^3 + 10 d^2 - 5 d - 36)$& $(\frac{5}{2}d^3+5d^2-\frac{5}{2}d-18)r-(\frac{25}{12}d^2+\frac{25}{4}d+2)((r-1)(r-2)-2+2n)$\\\hline
          $NP$ & $\frac{1}{2}(d^3+2d^2-d-4)$ & $\frac{1}{2}(d^3+2d^2-d-4)r-\frac{5}{12}d(d+3)((r-1)(r-2)-2+2n)$\\\hline
          $HF$ & $4 (11 d^3 - 68 d^2 + 79 d - 9)$ & $4 (11 d^3 - 68 d^2 + 79 d - 9)r-\frac{2}{3}(55 d^2 - 285 d + 132)((r-1)(r-2)-2+2n)$\\\hline
          $\text{BR}_N$ & $(d^3+2d^2-d-6)$ & $(d^3+2d^2-d-6)r-\frac{5}{6}d(d+3)((r-1)(r-2)-2+2n)$\\\hline
          $\text{BR}_T$ & $(13d^3-64d^2+53d+12)$ &$(13d^3-64d^2+53d+12)r-\frac{1}{6}(65 d^2 - 255 d + 36)((r-1)(r-2)-2+2n)$\\\hline
          $\text{BR}_F$ & $6 (13 d^3 - 79 d^2 + 83 d + 2)$ & $6 (13 d^3 - 79 d^2 + 83 d + 2)r-(65d^2-330d+111)((r-1)(r-2)-2+2n)$\\\hline
    \end{tabular}
    \caption{Enumerative numbers from intersecting the classes of moving curves $C_1$ and $C_2$ with divisor classes computed in Theorem A and Corollary 1.1}
    \label{tab:general_label}
\end{table}

\medskip\noindent
We may now intersect the curves $C_1,C_2$ with the divisor classes in Theorem A and Corollary 1.1. The following intersection numbers make the computations straightforward $$C_1 \cdot H=1, \quad C_1\cdot B=0,$$ and $$C_2\cdot H=r,\quad C_2\cdot B=2g(C)-2+2n = (r-1)(r-2)-2+2n.$$

\medskip\noindent 
Table \ref{tab:general_label} lists the intersection numbers of $C_1$ and $C_2$ with each of the divisors from the theorem and corollary, each of which is the answer to a distinct enumerative question.
Note, the conditions on $d$ to make the forgetful map dominant, guarantee that $d$ is divisible by $3$ which forces all of these values to be integers. 
\subsubsection{Example} Let $\Gamma\in \PP^{2[18]}$ generic and $C$ be a smooth curve of degree $5$ containing $\Gamma$. The moving curve $C_2$ is induced in $\PP^{2[18]}$ by letting $\Gamma$ vary in a general pencil in $|\mathcal{O}_C(\Gamma)|$. The number of times that 18 points in this pencil are the nodes of an irreducible curve of degree $9$ with a hyperflex is $C_2 \cdot HF=2252$.



\medskip\noindent



\section*{acknowledgments} \noindent
We would like to thank Izzet Coskun for useful conversations. During the preparation of this article the first author was partly supported by the CONACYT grant CB-2015/253061; he is currently a CONACYT Research Fellow in Mathematics, project No. 1036.


\begin{thebibliography}{MNOP06}
\bibitem[AC81]{AC} Arbarello, E. and Cornalba, M.: Footnotes to a paper of Beniamino Segre.
   \textit{Math. Ann} (256), 1981. 341--362.



\bibitem[BDPP13]{BDPP} Boucksom, S.; Demailly, J.P.; P?un, M.; Peternell, T.: The pseudo-effective cone of a compact Kähler manifold and varieties of negative Kodaira dimension. \textit{ Journal of Algebraic Geometry} (2013), 22 (2): 201?248.


\bibitem[DH88]{DH} Diaz, S. and Harris, J: Geometry of the {S}everi variety.
 \textit{Trans. Amer. Math. Soc.}(309),
    1988. 1--34.
     
     
\bibitem[DH88+]{DH2} Diaz, S. and Harris, J.: Ideals associated to deformations of singular plane curves.
\textit{Trans. Amer. Math. Soc.}(309) 1988. 433--468.
 

\bibitem[Fu98]{Fu98} Fulton, W.: Intersection theory. \textit{Springer-Verlag.} 1998.



\bibitem[Hu16]{Hu14} Huizenga, J.: Effective divisors on the Hilbert scheme of points in the plane and interpolation for stable bundles. \textit{J. Algebraic Geom.}(25) 2016. 19--75.


\bibitem[LR]{LR} Lozano Huerta, C. and Ryan, T: On the birational geometry of Hilbert schemes of points and Severi divisors. Preprint, 2019. https://arxiv.org/abs/1807.09881.



\bibitem[Se06]{Ser} Sernesi, E.: Deformations of algebraic schemes. Grundlehren der Mathematischen Wissenschaften (334), \textit{Springer-Verlag.} 2006.



\bibitem[Tr89]{Tre} Treger, R.: Plane curves with nodes. \textit{Can. J. Math.} Vol XLI, No. 2. (1989), 193--212.

\end{thebibliography}
\end{document}